\documentclass[12pt]{amsproc}
\usepackage{amsfonts,amssymb,latexsym, amsmath, amsthm, amsrefs}

\newtheorem{thm}{Theorem}[section]
\newtheorem{lem}[thm]{Lemma}
\newtheorem{cor}[thm]{Corollary}
\newtheorem{prop}[thm]{Proposition}

\newtheorem{conj}{Conjecture}

\theoremstyle{remark}

\newcommand{\GL}{{\mathrm {GL}}}
\newcommand{\rk}{{\mathrm {rk}}}
\newcommand{\SL}{{\mathrm {SL}}}

\newcommand{\PGL}{{\mathrm {PGL}}}

\def\N{{\mathbb N}}
\def\F{{\mathbb F}}

\def\leq{\leqslant}
\def\geq{\geqslant}

\def\Sp{{\rm Sp}}

\def\pr{{\mathrm{pr}}}
\def\H{{\/\mathrm H}}
\def\conj{{\mathrm c}}

\title{Metric ultraproducts of finite simple groups}

\author{Andreas Thom}
\address{A.T., Mathematisches Institut, U Leipzig,
PF 100920, 04009 Leipzig, Germany}
\email{andreas.thom@math.uni-leipzig.de}

\author{John S.Wilson}
\address{J.W., Mathematical Institute, Radcliffe Observatory Quarter,  Oxford OX2 6GG, England}
\email{john.wilson@maths.ox.ac.uk}

\begin{document}

\maketitle

\begin{abstract} Some new results on metric ultraproducts of finite simple groups are presented.  Suppose that $G$ is such a group, defined in terms of a non-principal ultrafilter $\omega$ on $\N$ and a sequence $(G_i)_{i\in\N}$ of finite simple groups, and that $G$ is neither finite nor a Chevalley group over an infinite field. Then $G$ is isomorphic to an ultraproduct of alternating groups or to an ultraproduct of finite simple classical groups.   The isomorphism type of $G$ determines which of these two cases arises, and, in the latter case, the $\omega$-limit of the characteristics of the groups $G_i$.  Moreover $G$ is a complete path-connected group with respect to the natural metric on $G$.
\bigskip

\noindent {\sc R\'esum\'e.} \  Nous pr\'esentons de nouveaux r\'esultats relatifs aux ultraproduits des groupes finis simples.  Soit $G$ un tel groupe, associ\'e \`a un ultrafiltre $\omega$ sur
$\N$ et une suite $(G_i)_{i\in\N}$ de groupes finis simples, et supposons que $G$ n'est ni fini ni un groupe de Chevalley sur un corps infini.  Un tel groupe $G$ est alors isomorphe soit \`a un ultraproduit de groups altern\'es, soit \`a un ultraproduit de groupes finis simples classiques.  La  classe d'isomorphisme de $G$ nous permet de distinguer ces deux cas, et dans le second cas, de d\'eterminer le $\omega$-limite des charact\'eristiques des groupes $G_i$.  Le groupe $G$ est de plus complet et connexe par arcs pour la m\'etrique naturelle sur $G$.
\end{abstract}

\section{Introduction}

Let $(G_i)_{i \in \N}$ be a sequence of groups and let $\omega$ be a non-principal ultrafilter on $\N$. The (algebraic) ultraproduct $\prod_{i \to \omega} G_i$ is defined to be the quotient of the Cartesian product $\prod_i G_i$ by the normal subgroup $\{ (g_i)_{i \in \N} \mid \{i \in \N \mid g_i=1\} \in \omega \}.$
We shall consider finite groups that are equipped with a normalized bi-invariant metric, that is, a metric $d \colon G \times G \to [0,1]$ such that $d(gh,gk)=d(h,k)=d(hg,kg)$ for all $g,h,k \in G$.
Given a sequence $(G_i,d_i)_{i \in \N}$ of such groups, we define the {\it metric} ultraproduct $\prod_{i \to \omega} (G_i,d_i)$ to be the quotient of $\prod_i G_i$ by the normal subgroup $$\left\{(g_i)_{i \in \N} \mid \lim_{i \to \omega} d_i(g_i,1)=0 \right\}.$$

Metric ultraproducts have been studied in a variety of situations, for example for tracial von Neumann algebras \cite{connes}, Banach spaces \cite{krivine} and more generally metric spaces; see \cite{wantiez} for a survey and more references. The study of metric ultraproducts of finite groups with bi-invariant metrics seems to be relatively recent \cite{elekszabo, thomas, stolzthom}.
Whereas the study of algebraic ultraproducts is closely related to the first-order theory of groups \cite{wilson}, the study of metric ultraproducts fits with the model theory for metric structures as developed in \cite{MR2436146}.

Every finite group with trivial centre has a natural bi-invariant metric
$$d_\conj(g,h) := \frac{\log |{\textrm{ccl}}(gh^{-1})|}{\log |G|},$$
where ${\textrm{ccl}}(g)\!:= \{kgk^{-1} \mid k \in G\}$. We call $d_\conj$ the {\em conjugacy metric} on $G$. Other natural examples of bi-invariant metrics arise from combinatorics and geometry. For the symmetric group $S_n$ the {\em Hamming distance} 
$d_\H \colon S_n \times S_n \to [0,1]$ is defined by
$d_\H(\sigma,\tau):= \frac1n |\{1 \leq i \leq n \mid \sigma(i) \neq \tau(i) \}|$.  For general linear groups, the {\em projective rank distance} $d_\pr \colon \GL_n(F) \times \GL_n(F) \to [0,1]$ defined by
$d_\pr(g,h):= \frac1n \left( \min_{\alpha\in F^\times} {\rm rk}(g-\alpha h)\right)$ is a pseudometric and it induces a metric on $\PGL_n(F)$ and its subgroups.  Thus we have natural metrics on the alternating groups and all finite simple classical groups, and they turn out to be asymptotically equivalent (as $n\to_\omega\infty$) to the conjugacy metric \cite{stolzthom}. The proper  normal subgroups in any algebraic ultraproduct of finite simple groups are totally ordered by inclusion and have a maximal element; the quotient by this element is precisely the metric ultraproduct and so the metric ultraproduct is a simple group \cite{stolzthom}. 

Metric ultraproducts of symmetric groups of unbounded degree with respect to the Hamming distance have gained attention recently, since a countable group is sofic \cite{gromov, weiss} if and only it can be embedded such groups; see \cite{elekszabo}.

\section{The isomorphism problem for metric ultraproduct of finite simple groups}

Let $G$ be a metric ultraproduct of a sequence of non-abelian finite simple groups.
Using the classification of non-abelian finite simple groups and a result of Point \cite{Point}, one can prove that $G$ will be either finite, a simple group of Lie type over an ultraproduct of finite fields, a metric ultraproduct of alternating groups $(A_{n_i})_{i \in \N}$ with $\lim_{i \to \omega} n_i = \infty$, or a metric ultraproduct of finite simple classical groups of unbounded ranks.  We recall that the finite simple classical groups are projective special linear, unitary, symplectic and orthogonal groups and they are Chevalley groups $A_n(q), B_n(q), C_n(q), D_n(q), {}^2A_n(q^2),$ ${}^2D_n(q^2)$, where $n$ is the rank and $q$ is a power of a prime $p$, the characteristic of the group. 

For any  sequence $(G_i)_{i \in \N}$ of finite simple classical groups with characteristics $p_i$, 
we define the characteristic  of the sequence 
(with respect to a given ultrafilter $\omega$ on $\N$) to be  
$\lim_{i \to \omega} p_i \in {\mathbb P} \cup \{\infty\}$. Here, $\mathbb P$ denotes the set of prime numbers.

There is no difference between the metric ultraproducts of special linear groups and general linear groups.
\begin{lem} \label{dense}
Let $G$ be a metric ultraproduct of general linear groups of unbounded rank.  Then $G$ is equal to the corresponding metric ultraproduct of special linear groups.  Thus $G$ is perfect and indeed every element can be written as a commutator.
\end{lem}

Our main result is:

\begin{thm} \label{main}
Let $(G_i)_{i \in \N}$ be a sequence of non-abelian finite simple groups. The isomorphism type of\/ $\prod_{i \to \omega} (G_i,d_\conj)$ determines whether it is isomorphic to a metric ultraproduct of alternating groups or a metric ultraproduct of simple classical groups. In the latter case, the isomorphism type also determines the characteristic of the sequence $(G_i)_{i \in \N}$.
\end{thm}

As a very special case we have

\begin{cor}
The groups $\prod_{i \to \omega} (A_n,d_\H)$ and all groups $\prod_{i \to \omega} A_n(p)$ with $p$ prime are pairwise non-isomorphic.
\end{cor}

The proof relies on the study of the structure of centralizers of elements of prime order $p$ in the metric ultraproduct $G$ of a sequence $(G_i)_{i \in \N}$. Let $\chi$ be the characteristic of this sequence.   
 Roughly speaking, if $p\neq \chi$ then the centralizer of  each element of order $p$ in $G$ can be shown to have a normal subgroup which is a product of quasisimple groups (whose simple quotients are themselves metric
ultraproducts of finite simple classical groups). On the other hand, if $p=\chi$, then there is a specific element of order $p$ whose centralizer is essentially a semi-direct product of an elementary abelian $p$-group by a product of simple groups; and every infinite normal subgroup contains this elementary abelian normal subgroup. In the case of alternating groups, the second type of behaviour can occur for all primes $p$. These facts allow us to distinguish the case of alternating groups from the case of classical groups, and to determine the characteristic of the sequence $(G_i)_{i \in \N}$ from the isomorphism type of the metric ultraproduct.

We outline the arguments in the cases of alternating groups and projective special linear groups.  (The arguments for other classical groups are more complicated.)  First, we use the following result to show that an element of finite order $k$ in a metric ultraproduct of projective general linear groups of characteristic prime to $k$ can be represented by a sequence of elements which are essentially of order $k$.

\begin{prop} \label{preparation}
Let $k, n \in \N$, let $F$ be a field with $({\rm char}\,F,k)=1$ and let $V$ be the natural module for $\GL_n(F)$. 
Let $y \in \GL_n(F)$ and $\alpha \in F^{\times}$. There exist $x \in \GL_n(F)$ and $x$-invariant subspaces $L$, $S$ of $V$ such that
$$V=L\oplus S,\quad (x|_L)^k=\alpha\cdot 1_L \quad\hbox{and}\quad (x|_S)=1_S\eqno(1)$$ and
$$\max\{\dim S, \rk(x-y) \} \leq \rk(y^k-\alpha 1_n).$$
\end{prop}

Secondly, we need to know that any element of the centralizer in the metric ultraproduct lies in the image of the product of centralizers inside the Cartesian product.

\begin{prop} \label{centralizer} Let $k, n \in \N$, let $F$ be a field with $({\rm char}\,F,k)=1$ and let $x \in \GL_n(F)$.
Suppose that the natural module $V$ has $x$-invariant subspaces $L$, $S$ satisfying {\rm(1)} for some
$\alpha\in F_q^{\times}$.  
Let $\phi \in \GL_n(F)$. 
Then
there exists $\psi \in \GL_n(F)$ with 
$$x \psi =\psi x\quad \mbox{and} \quad
\rk(\phi - \psi) \leq 2k^2 \cdot \rk (x \phi - \phi x)+3\dim S.$$
\end{prop}

Results of the following type are  well known.

\begin{prop} \label{general} Let $k, n \in \N$,and let $F$ be a field with $({\rm char}\,F,k)=1$. Let $\alpha\in F^\times$ and $x \in \GL_n(q)$.  If the natural module $V$ for $\GL_n(q)$ has $x$-invariant subspaces $L$, $S$ satisfying {\rm(1)}, then the centralizer of $x$ is a product of at most $k+1$ general linear groups.
\end{prop}

On the other hand, if $p$ is equal to the characteristic of the field, the situation is different.

\begin{prop}\label{niceblock} Let $p$ be a prime and $F$ a field of characteristic $p$, and let $n\geq2$.  Let $G$ be either of the groups 
$\SL_{2n}(F)$, $\Sp_{2n}(F)$.  Then $G$ has an element $x$ of order $p$ with $d_\pr(x,1_{2n})=\frac12$ whose centralizer  is a semidirect product $A\rtimes H$ having the following properties$:$
\begin{enumerate}\item[\rm(i)]  $A$ is an elementary abelian $p$-group and $H$ is a quasisimple classical group$;$
\item[\rm(ii)] both $A$, $H$ contain elements of projective rank length at least $\frac12;$
\item[\rm(iii)] the set of commutators $[u,h]=u^{-1}h^{-1}uh$ with $u\in A$, $h\in H$ contains elements of projective rank length at least $\frac13(1-2 n^{-1})$.  \end{enumerate} \end{prop}

This leads to the following theorem.

\begin{thm} Let $\omega$ be a non-principal ultrafilter on $\N$.
Let $(F_i)$  be a sequence of fields, let $p_i={\mathrm {char}}\, F_i$ for each $i$ and let $G=\prod_{i \to \omega} A_{n_i}(F_i)$.  Write $$\chi=\displaystyle \lim_{i \to \omega} p_i.$$  
\begin{enumerate} \item[\rm(a)]  If $x\in G$ has prime order $p\neq \chi$ then its centralizer ${\rm C}_G(x)$ is an extension of its centre by a direct product of at most $p$ simple groups.
\item[\rm(b)]  If $\chi$ is finite then $G$ has an element $y$ of order $\chi$ such that ${\rm C}_G(y)$ has a non-trivial normal elementary abelian group on which ${\rm C}_G(y)$ acts as an infinite group. \end{enumerate} 

Consequently, the isomorphism class of $G$ determines $\lim_{i \to \omega} p_i \in {\mathbb P} \cup \{\infty\}.$
\end{thm}

The following result for ultraproducts of alternating groups finishes the proof of Theorem \ref{main}.

\begin{thm}  Let $G=\prod_{i \to \omega} (A_{n_i},d_{\rm H})$ be a metric ultraproduct of alternating groups
with $n_i\to_\omega \infty$.  Let $x$ be an element of $G$ of prime order $p$ and let $C={\rm C}_G(x)$.  Then $C$ has the form $(M\rtimes T_1)\times T_2$ where  $T_1$ is a metric
ultraproduct of alternating groups of unbounded degree, where $M$ is an infinite $\F_pT_1$-module which is simple modulo a submodule of order $p$ and where $T_2$ is a metric
ultraproduct of alternating groups of unbounded degree if $\ell_{\rm H}(x)<1$ and $T_2$ is finite if
$\ell_{\rm H}(x)=1$. 
\end{thm} 

It remains to decide whether metric ultraproducts of classical groups of different types can be isomorphic.

\section{Geometry of metric ultraproducts}

Let $(G_i,d_i)_{i \in \N}$ be a sequence of groups with normalized bi-invariant metrics. Then the metric ultraproduct $G=\prod_{i \to \omega} (G_i,d_i)$ is again a metric space, with metric induced by the pseudometric 
$$d((g_i),(h_i)):= \lim_{i \to \omega} d_i(g_i,h_i)$$
on $\prod_{i\to\omega}G_i$.
It is natural to study topological and metric properties of $G$.  The metric geometry of $G$ reflects asymptotic features of the behaviour of the metrics $d_i$.
When the groups $G_i$ are finite simple groups, the most natural choices of bi-invariant metrics are asymptotically equivalent and so the topology on $G$, unlike the metric, is independent of the choices.  

\begin{thm}\label{connected}
Let $(G_i)_{i \in \N}$ be a sequence of finite simple groups. Then the metric ultraproduct is either discrete, in which case it is finite or a classical group of Lie type over an ultraproduct of finite fields, or complete and path-connected.
\end{thm}

The non-trivial cases of the preceding theorem again reduce to the study of metric ultraproducts of sequences in one of the families. 
Completeness is a consequence of more general facts; see \cite[p.\ 147]{wantiez}.

For metric ultraproducts of alternating groups and projective special linear groups we can prove a stronger assertion than Theorem \ref{connected}.
We recall that a metric space $(X,d)$ is called {\em geodesic} if, for each pair of points $x,y \in X$, there is an isometric embedding $\iota \colon [0,d(x,y)] \to X$ such that $\iota(0)=x$ and $\iota(d(x,y))=y$.

\begin{thm}
The non-discrete metric ultraproducts of sequences $(A_{n_i},d_\H)_{i \in \N}$ and $(A_{n_i}(q_i),d_\pr)_{i \in \N}$ for prime powers $q_i$ are complete geodesic metric spaces.
\end{thm}

\section*{Acknowledgments}

The first author was supported by the ERC. This project was started during a visit of the second author to the Universit\"at Leipzig as guest of the Felix-Klein-Colleg.



\end{document}